\documentclass[12pt]{article}
\usepackage{amsmath,amssymb,amsthm,amscd,latexsym}
\renewcommand{\.}{\hbox{.}~}
\renewcommand{\le}{\leqslant}
\renewcommand{\ge}{\geqslant}
\newcommand{\wt}{\widetilde}
\renewcommand{\[}{\lfloor}
\renewcommand{\]}{\rfloor}
\renewcommand{\d}{\mathrm{d}}

\DeclareMathOperator{\ord}{ord}
\DeclareMathOperator{\Real}{Re}
\renewcommand{\Re}{\Real}

\providecommand{\doublesb}[2]{_{\genfrac{}{}{0pt}{1}{#1}{#2}}}
\newcommand{\ba}{\boldsymbol a}
\newcommand{\bn}{\boldsymbol n}
\newcommand{\bl}{\boldsymbol l}
\begingroup
\catcode`\@=11
\global\tagsleft@false
\endgroup
\begin{document}
\theoremstyle{plain}
\newtheorem{theorem}{Theorem}
\newtheorem{prop}{Proposition}
\theoremstyle{remark}
\newtheorem*{remark}{Remark}
\newtheorem*{acknowledgements}{Acknowledgements}
\title{An essay on irrationality measures\\of $\pi$ and other logarithms}
\author{Wadim Zudilin%
\thanks{The work is supported by an Alexander von Humboldt
research fellowship and partially supported by grant
no.~03-01-00359 of the Russian Foundation for Basic Research.}
\quad\rm(Moscow)}
\date{%
\hbox to\hsize{\vbox{%
\centerline{E-print \texttt{math.NT/0404523}}
\smallskip
\centerline{14 May 2004}
}}}
\maketitle
\rightline{\it To my teacher and friend A.\,I.~Galochkin}
\rightline{\it on the occasion of his 60th birthday}
\bigskip
\pagestyle{myheadings}
\markboth{W.~Zudilin}{Irrationality measures of logarithms}
\thispagestyle{empty}

Let $a\in\mathbb Q\cap(0,2]$, $a\ne1$. Then the sequence
of quantities
\begin{equation}
\int_0^1\frac{x^n(1-x)^n}{(1-(1-a)x)^{n+1}}\,\d x
\in\mathbb Q\log a+\mathbb Q,
\qquad n=0,1,2,\dots,
\label{eq:1}
\end{equation}
produces `good' rational approximations to $\log a$.
There are several ways of performing integration in~\eqref{eq:1}
in order to show that the integral lies in $\mathbb Q\log a+\mathbb Q$;
we give an exposition of different methods below.
The aim of this essay is to demonstrate how suitable
generalizations of the integrals in~\eqref{eq:1} allow
to prove the best known results on irrationality measures
of the numbers $\log2$, $\pi$ and $\log3$. Although methods
presented below work in general situations (e.g., for
certain $\mathbb Q$-linear forms in logarithms) as well, the three
numbers seem to be very nice and important models
for our exposition.

Bounds for irrationality measures are presented
by means of upper estimates for irrationality exponents.
Recall that the {\it irrationality exponent\/}
of a real irrational number~$\gamma$ is defined
by the relation
\begin{align*}
\mu=\mu(\gamma)
=\inf
&\{c\in\mathbb R:\text{the inequality $|\gamma-a/b|\le|b|^{-c}$ has}
\\ &\qquad
\text{only finitely many solutions in $a,b\in\mathbb Z$}\}.
\end{align*}
The estimates for $\mu(\gamma)$ are deduced by constructing
sequences of linear forms involving $\gamma$ and using
standard tools of the following shapes.

\begin{prop}[\rm\cite{Ha1}, Lemma 3.1]
\label{pr:1}
Let $\gamma\in\mathbb R$ be irrational.
Suppose that a sequence of linear forms $b_nx-a_n$,
with integer coefficients from the field of rationals
or an imaginary quadratic field, satisfies
$$
\limsup_{n\to\infty}\frac{\log|b_n|}n\le C_1,
\qquad
\lim_{n\to\infty}\frac{\log|b_n\gamma-a_n|}n=-C_0
$$
for some positive real $C_0$ and $C_1$. Then
$\mu(\gamma)\le1+C_1/C_0$.
\end{prop}

\begin{prop}[\rm\cite{Ha2}, Lemma 2.1]
\label{pr:2}
Let $\omega,\omega'\in\mathbb R$ be two irrational numbers.
Suppose that sequences of linear forms $b_nx-a_n$ and $b_nx-a_n'$,
with integer coefficients from the field of rationals
or an imaginary quadratic field, satisfies
$$
\limsup_{n\to\infty}\frac{\log|b_n|}n\le C_1,
\qquad
\lim_{n\to\infty}\frac{\log|b_n\omega-a_n|}n=-C_0,
\quad
\lim_{n\to\infty}\frac{\log|b_n\omega'-a_n'|}n=-C_0'
$$
for some positive real constants $C_0<C_0'$ and $C_1$.
Then any nonzero element $\gamma\in\mathbb Q\omega+\mathbb Q\omega'$
is irrational with the bound
$\mu(\gamma)\le1+C_1/C_0$ for the irrationality exponent.
\end{prop}

\begin{remark}
In fact, the statement of Lemma~2.1 in \cite{Ha2} slightly
differs from our last claim, but one can easily verify
that the proof given there proves our `modification' as well.
\end{remark}

\section{Irrationality measure for $\log2$ (after E.~Rukhadze)}
\label{sec:1}

\subsection{Gauss hypergeometric function}
\label{sec:1.1}
It is worth performing a slightly general integral than~\eqref{eq:1},
namely
\begin{equation}
I(m,n_0,n_1;a)
=\int_0^1\frac{x^{n_0}(1-x)^{n_1}}{(1-(1-a)x)^{m+1}}\,\d x
\label{eq:2}
\end{equation}
for non-negative integers $m,n_0,n_1$,
provided the condition $\max\{m,n_0\}\le n_1$ holds for further
convenience. The integral in~\eqref{eq:2} is exactly
Euler's integral for the Gauss hypergeometric series:
\begin{align}
I(m,n_0,n_1;a)
&=\frac{\Gamma(n_0+1)\,\Gamma(n_1+1)}{\Gamma(n_0+n_1+2)}
{}_2F_1\biggl(\begin{matrix} m+1, \, n_0+1 \\ n_0+n_1+2 \end{matrix}
\biggm|1-a\biggr)
\nonumber\\
&=\frac{\Gamma(n_1+1)}{\Gamma(m+1)}
\sum_{\nu=0}^\infty\frac{\Gamma(m+1+\nu)\,\Gamma(n_0+1+\nu)}
{\Gamma(1+\nu)\,\Gamma(n_0+n_1+2+\nu)}(1-a)^\nu
\label{eq:3}
\end{align}
(see, e.g., \cite{AAR}, Section~2.2). The latter sum may be written
as
\begin{equation}
I(m,n_0,n_1;a)
=\sum_{\nu=0}^\infty R(\nu)(1-a)^\nu,
\label{eq:4}
\end{equation}
where
\begin{equation}
R(t)=\frac{(t+1)(t+2)\dotsb(t+m)}{m!}
\cdot\frac{n_1!}{(t+n_0+1)(t+n_0+2)\dotsb(t+n_0+n_1+1)}
\label{eq:5}
\end{equation}
and $R(t)=O(t^{-1})$ as $t\to\infty$ by $m\le n_1$.
Denote $m^*=\min\{m,n_0\}$ and $n_0^*=\max\{m,n_0\}$
and decompose the rational function~\eqref{eq:5} in a sum
of partial fractions:
\begin{equation}
R(t)=\sum_{k=n_0^*}^{n_0+n_1}\frac{A_k}{t+k+1}
=\sum_{k=n_0^*}^{n_0+n_1}
\frac{(-1)^{m+n_0-k}\binom km\binom{n_1}{k-n_0}}{t+k+1}.
\label{eq:6}
\end{equation}
Then by~\eqref{eq:4} we obtain
\begin{align}
I(m,n_0,n_1;a)
&=\sum_{\nu=-m^*}^\infty R(\nu)(1-a)^\nu
=\sum_{k=n_0^*}^{n_0+n_1}A_k(1-a)^{-(k+1)}
\sum_{\nu=-m^*}^\infty\frac{(1-a)^{\nu+k+1}}{\nu+k+1}
\nonumber\\
&=\sum_{k=n_0^*}^{n_0+n_1}A_k(1-a)^{-(k+1)}
\biggl(\sum_{l=1}^\infty-\sum_{l=1}^{k-m^*}\biggr)\frac{(1-a)^l}l
\nonumber\\
&=-\log a\cdot\sum_{k=n_0^*}^{n_0+n_1}A_k(1-a)^{-(k+1)}
-\sum_{k=n_0^*}^{n_0+n_1}\sum_{l=1}^{k-m^*}\frac{A_k(1-a)^{l-(k+1)}}l,
\label{eq:7}
\end{align}
hence
\begin{equation}
I(m,n_0,n_1;a)(1-a)^{n_0+n_1+1}\cdot d^{n_0+n_1-m^*}D_{n_0+n_1-m^*}
\in\mathbb Z\log a+\mathbb Z,
\label{eq:8}
\end{equation}
where $d$ denotes the denominator of~$a$ and
$D_n$ stands for the least common multiple of the numbers
$1,2,\dots,n$. By the prime number theorem, we have
the following asymptotic formula:
$$
\lim_{n\to\infty}\frac{\log D_n}n=1.
$$

\subsection{Arithmetic valuation}
\label{sec:1.2}
The inclusion \eqref{eq:8} may be essentially improved
in several cases, and it is the observation that allowed
Rukhadze to prove the record irrationality measure for~$\log2$.

The symmetry of the ${}_2F_1$-series in~\eqref{eq:3}
with respect to its upper parameters $m+1$ and $n_0+1$
gives us a way to write the identity
\begin{equation}
\frac{I(m,n_0,n_1;a)}{\Gamma(n_0+1)\,\Gamma(n_1+1)}
=\frac{I(n_0,m,n_0+n_1-m;a)}{\Gamma(m+1)\,\Gamma(n_0+n_1-m+1)}
\label{eq:9}
\end{equation}
(which is not so evident if one looks on definition~\eqref{eq:2}).
The inclusion~\eqref{eq:8} written for the $I$-quantity on the right
of~\eqref{eq:9},
$$
I(n_0,m,n_0+n_1-m;a)(1-a)^{n_0+n_1+1}\cdot d^{n_0+n_1-m^*}D_{n_0+n_1-m^*}
\in\mathbb Z\log a+\mathbb Z,
$$
and the equality
\begin{align*}
&
I(m,n_0,n_1;a)(1-a)^{n_0+n_1+1}\cdot d^{n_0+n_1-m^*}D_{n_0+n_1-m^*}
\cdot\frac{m!\,(n_0+n_1-m)!}{n_0!\,n_1!}
\\ &\qquad
=I(n_0,m,n_0+n_1-m;a)(1-a)^{n_0+n_1+1}\cdot d^{n_0+n_1-m^*}D_{n_0+n_1-m^*}
\end{align*}
imply that if $\Phi(m,n_0,n_1)$ is the denominator
of the quotient
$$
\frac{m!\,(n_0+n_1-m)!}{n_0!\,n_1!},
$$
then
\begin{equation}
I(m,n_0,n_1;a)(1-a)^{n_0+n_1+1}\cdot d^{n_0+n_1-m^*}D_{n_0+n_1-m^*}
\cdot\Phi(m,n_0,n_1)^{-1}
\in\mathbb Z\log a+\mathbb Z.
\label{eq:10}
\end{equation}

By the well-known formula, for each prime~$p$ we have
$\ord_pN!=\[N/p\]+\[N/p^2\]+\[N/p^3\]+\dotsb$,
where $\[\,\cdot\,\]$
denotes the integral part of a number. Therefore
\begin{equation}
\Phi(m,n_0,n_1)
=\prod_pp^{\phi(p)+\phi(p^2)+\phi(p^3)+\dotsb},
\label{eq:11}
\end{equation}
where
$$
\phi(t)
=\max\biggl\{0,\biggl\[\frac{n_0}t\biggr\]+\biggl\[\frac{n_1}t\biggr\]
-\biggl\[\frac mt\biggr\]-\biggl\[\frac{n_0+n_1-m}t\biggr\]\biggr\}.
$$
The final remark (made by G.~Chudnovsky in~\cite{Ch} together
with introducing the method of asymptotic evaluation of
the factors like~\eqref{eq:11}) consists in the fact that
the divisor
\begin{equation}
\wt\Phi(m,n_0,n_1)
=\prod_{p>\sqrt{n_1}}p^{\phi(p)}
\label{eq:12}
\end{equation}
of $\Phi(m,n_0,n_1)$ gives the main contribution in the asymptotic
of~\eqref{eq:11} and may be easily controlled.

\subsection{Irrationality result}
\label{sec:1.3}
The choice $a=2$ and $n_0=6n$, $m=7n$, $n_1=8n$, where $n$ is the positive
integer parameter increasing to $\infty$, allowed E.~Rukhadze in~\cite{Ru}
to prove the following result (see also \cite{Ha1}, \cite{Vi} and \cite{Br}).

\begin{theorem}
\label{th:1}
The irrationality exponent of $\log2$ satisfies the inequality
$$
\mu(\log2)\le3.89139977\dotsc.
$$
\end{theorem}

We will briefly indicate required ingredients of the proof.
For the above choice of the parameters we set
$$
I_n=I(7n,6n,8n;2)
=\int_0^1\biggl(\frac{x^6(1-x)^8}{(1+x)^7}\biggr)^n\frac{\d x}{1+x}
=\bar A_n\log2-\bar B_n,
$$
where, by \eqref{eq:6} and \eqref{eq:7},
$$
\bar A_n=(-1)^n\sum_{k=7n}^{14n}\binom k{7n}\binom{8n}{k-6n}.
$$
Then
\begin{align}
\lim_{n\to\infty}\frac{\log I_n}n
&=\log\max_{0<x<1}\frac{x^6(1-x)^8}{(1+x)^7}
\nonumber\\
&=\log\frac{2^53^3(7734633\sqrt{393}-153333125)}{7^7}
=-11.84497806\dots
\label{eq:13}
\end{align}
and, thanks to Stirling's asymptotic formula for the factorial,
\begin{align}
\lim_{n\to\infty}\frac{\log|\bar A_n|}n
&=\lim_{n\to\infty}\frac1n\log\max_{7n\le k\le14n}
\binom k{7n}\binom{8n}{k-6n}
\nonumber\\
&=\log\max_{7<y<14}\biggl(\frac{y^y}{7^7(y-7)^{y-7}}
\cdot\frac{8^8}{(y-6)^{y-6}(14-y)^{14-y}}\biggr)
\nonumber\\
&=\log\frac{2^53^3(7734633\sqrt{393}+153333125)}{7^7}
=12.68147230\dotsc.
\label{eq:14}
\end{align}
Concerning the asymptotic behaviour of the value
$\Phi_n=\wt\Phi(7n,6n,8n)$ in \eqref{eq:12}, we use
the fact $\phi(t)=\varpi_0(n/t)$, where
\begin{align*}
\varpi_0(x)
&=\max\bigl\{0,\[6x\]+\[8x\]-2\[7x\]\bigr\}
\\
&=\begin{cases}
1 & \text{if $x\in\bigl[\frac18,\frac17\bigr)
\cup\bigl[\frac14,\frac27\bigr)
\cup\bigl[\frac38,\frac37\bigr)
\cup\bigl[\frac12,\frac47\bigr)
\cup\bigl[\frac23,\frac57\bigr)
\cup\bigl[\frac56,\frac67\bigr)$}, \\
0 & \text{otherwise}.
\end{cases}
\end{align*}
Therefore,
\begin{equation}
\lim_{n\to\infty}\frac{\log\Phi_n}n
=\int_0^1\varpi_0(x)\d\psi(x)
=\log\frac{2^{15}3^3}{7^7}+\frac{\pi(3+6\sqrt2-4\sqrt3)}6
=2.45775406\dots,
\label{eq:15}
\end{equation}
where $\psi(x)$ denotes the logarithmic derivative
of the gamma function. Using inclusions~\eqref{eq:10}
and the asymptotics \eqref{eq:13}--\eqref{eq:15}, we obtain
\begin{align*}
C_0&=-\log(7734633\sqrt{393}-153333125)
+10\log2-8+\frac{\pi(3+6\sqrt2-4\sqrt3)}6
\\
&=6.30273213\dots,
\\
C_1&=\log(7734633\sqrt{393}+153333125)
-10\log2+8-\frac{\pi(3+6\sqrt2-4\sqrt3)}6
\\
&=18.22371823\dots,
\end{align*}
in the notation of Proposition~\ref{pr:1} and, finally,
conclude with the estimate
$$
\mu(\log2)\le1+\frac{C_1}{C_0}
=3.89139977\dotsc.
$$

The result for the measure of $\log2$
may be compared with that obtained in
simpler settings $n_0=n_1=m=n$ (as in~\eqref{eq:1}):
$$
C_0=-2\log(\sqrt2-1)-1=2\log(\sqrt2+1)-1,
\qquad
C_1=2\log(\sqrt2+1)+1,
$$
hence
$$
\mu(\log2)\le1+\frac{C_1}{C_0}
\le1+\frac{2\log(\sqrt2+1)+1}{2\log(\sqrt2+1)-1}
=4.62210083\dotsc.
$$

\section{Irrationality measure for $\pi$ (after M.~Hata)}
\label{sec:2}

\subsection{Simultaneous approximations to logarithms}
\label{sec:2.1}
The change of variable $z=1-(1-a)x$
in~\eqref{eq:1} transforms the integral~\eqref{eq:1} into
\begin{equation}
\frac{(-1)^{n+1}}{(1-a)^{2n+1}}
\int_1^a\frac{(z-1)^n(z-a)^n}{z^{n+1}}\,\d z.
\label{eq:16}
\end{equation}
Instead of decomposing the latter integral we will perform
a more general {\it complex\/} integral
$$
I_k(\ba,m,\bn;a)=\int_{\Gamma_{1,a}}
\frac{(z-1)^{n_0}(z-a_1)^{n_1}\dotsb(z-a_k)^{n_k}}{z^{m+1}}\,\d z,
$$
where $\Gamma_{1,a}$ denotes a smooth oriented path from~$1$ to~$a$
contained in $\mathbb C\setminus\{0\}$; the parameters $a,a_1,\dots,a_k$
are complex numbers distinct from $0,1$; the exponents
$n_0,n_1,\dots,n_k,m$
are positive integers. The integral in~\eqref{eq:16}
corresponds to $k=1$, $a_1=a$ and $n_0=n_1=m=n$.
Setting additionally $a_0=1$, we may compute, as in~\cite{Ha2}, Section~3,
\begin{align}
I_k(\ba,m,\bn;a)
&=\sum_{l_0=0}^{n_0}\sum_{l_1=0}^{n_1}\dotsi\sum_{l_k=0}^{n_k}
A_{\bl}\binom{n_0}{l_0}\binom{n_1}{l_1}\dotsb\binom{n_k}{l_k}
\int_{\Gamma_{1,a}}z^{l_0+l_1+\dots+l_k-m-1}\,\d z
\nonumber\\
&=\mathop{\sum\dotsi\sum}_{l_0+\dots+l_k\ne m}
\frac{A_{\bl}}{l_0+\dots+l_k-m}
\binom{n_0}{l_0}\dotsb\binom{n_k}{l_k}
(a^{l_0+\dots+l_k-m}-1)
\nonumber\\ &\qquad
+\mathop{\sum\dotsi\sum}_{l_0+\dots+l_k=m}
A_{\bl}\binom{n_0}{l_0}\dotsb\binom{n_k}{l_k}\cdot\log a,
\label{eq:17}
\end{align}
where
$$
A_{\bl}=A_{l_0,l_1,\dots,l_k}
=(-1)^{l_0+l_1+\dots+l_k}a_1^{n_1-l_1}\dotsb a_k^{n_k-l_k}
$$
and we use the formula
$$
\int_{\Gamma_{1,a}}z^{l-1}\,\d z
=\int_1^az^{l-1}\,\d z
=\begin{cases}
a^l/l & \text{if $l\ne0$}, \\
\log a & \text{if $l=0$}.
\end{cases}
$$

The main idea is that the coefficient of $\log a$
in the linear form~\eqref{eq:17} does not depend on the choice of~$a$
(but of course the analytic behaviour of the integral does!).
The suitable and natural choice of~$a$
is from the set $\{a_1,\dots,a_k\}$. Then the above
quantities $I_k$ produce simultaneous approximations
to $\log a_1,\dots,\log a_k$.

\subsection{Analytic and arithmetic ingredients}
\label{sec:2.2}
Our basic consideration will be devoted to the case $k=2$,
which is used in~\cite{Ha2} to give the linear independence
measure of $\pi$ and $\log2$ over~$\mathbb Q$ (in particular,
the irrationality measure of~$\pi$)
and the new irrationality measure of~$\pi/\sqrt3$.

Thus, Hata~\cite{Ha2} takes $k=2$ (that really
gives an extension of~\eqref{eq:16},
and hence of~\eqref{eq:1}) and substitute $a=a_1$ and $a=a_2$ to get
nice simultaneous approximations to $\log a_1$ and $\log a_2$.
Hata `restricts' himself from the beginning to
considering the particular case $n_0=n_1=n_2=2n$ and $m=3n$,
where $n$ is an increasing parameter. However, this simple choice
produces the best possible number-theoretic results,
and our consideration of the general case
$$
n_0=\alpha_0n, \quad
n_1=\alpha_1n, \quad
n_2=\alpha_2n, \quad
m=\alpha n,
$$
where $\alpha_0,\alpha_1,\alpha_2,\alpha$ are positive integers,
is mostly due to methodological reasons.

Write the integrals in the form
\begin{equation}
J_{j,n}=I_2(a_j)=\int_{\gamma_j}\frac{e^{nf(z)}}z\,\d z,
\qquad j=1,2,
\label{eq:18}
\end{equation}
where
$$
f(z)=\alpha_0\log(z-a_0)+\alpha_1\log(z-a_1)
+\alpha_2\log(z-a_2)-\alpha\log z
$$
and the path $\gamma_j$ joints the points $1$ and $a_j$
and goes through the corresponding saddle point.
The saddle points $\xi_0,\xi_1,\xi_2$
are solutions of the equation $f'(z)=0$
becoming the cubic polynomial equation: two of these
saddles correspond to the growth of the integrals in~\eqref{eq:18},
$$
\lim_{n\to\infty}\frac{\log|J_{1,n}|}n=\Re f(\xi_1),
\qquad
\lim_{n\to\infty}\frac{\log|J_{2,n}|}n=\Re f(\xi_2),
$$
while the third saddle $\xi_0$ determines the asymptotic
behaviour of the coefficients of the linear forms.

To compute the arithmetic of the coefficients
we should evaluate the true denominators of the products
$$
\frac1{l_0+l_1+l_2-\alpha n}
\binom{\alpha_0n}{l_0}\binom{\alpha_1n}{l_1}\binom{\alpha_2n}{l_2},
\qquad l_0+l_1+l_2\ne\alpha n.
$$
Clearly the least common multiple $D_{\beta n}$,
where $\beta=\max\{\alpha,\alpha_0+\alpha_1+\alpha_2-\alpha\}$,
is required but some primes $p>\sqrt{Cn}$ may be then excluded
from this $D_{\beta n}$ by considering the following problem:
determine primes $p$ dividing all the integers
$$
\binom{\alpha_0n}{l_0}\binom{\alpha_1n}{l_1}\binom{\alpha_2n}{l_2}
$$
under the additional condition $l_0+l_1+l_2\equiv\alpha n\pmod p$.
Writing $x=\{n/p\}$ and $y_j=\{l_j/p\}$, $j=0,1,2$, for the
fractional parts, we reduce the problem to minimizing
the $1$-periodic integer-valued function
$$
\varpi(x,y_0,y_1,y_2)
=\sum_{j=0}^2\bigl(\[\alpha_jx\]-\[y_j\]-\[\alpha_jx-y_j\]\bigr)
$$
on the cube $(y_0,y_1,y_2)\in[0,1)^3$ under the additional
hypothesis $y_0+y_1+y_2\equiv\alpha x\pmod1$.
(The last condition means that knowledge of $x,y_0,y_1$
determines the remaining value $y_2$ uniquely.) Denote by $\varpi_0(x)$
the required minimum. For example, Hata's choice
$\alpha_0=\alpha_1=\alpha_2=2$, $\alpha=3$ gives
$$
\varpi_0(x)=\begin{cases}
1 & \text{if $x\in\bigl[\frac12,\frac23\bigr)$}, \\
0 & \text{otherwise}.
\end{cases}
$$

There is also a `problem' of finding the true denominators of
$A_{\bl}$ and $A_{\bl}a^{l_0+l_1+l_2-m}$.
For example, in the case $a_1=2$, $a_2=1+i$ (of simultaneous
approximations to $\log2$ and $\pi$) we have
\begin{align*}
&
(-1)^{l_0+l_1+l_2}A_{\bl}a_0^{l_0+l_1+l_2-m}
=2^{n_1-l_1}(1+i)^{n_2-l_2}
\in\mathbb Z[i],
\\ &
(-1)^{l_0+l_1+l_2}A_{\bl}a_1^{l_0+l_1+l_2-m}
=2^{n_1+l_0+l_2-m}(1+i)^{n_2-l_2}
\\ &\qquad
=2^{n_1+l_0-m}(1+i)^{l_2}(1-i)^{l_2}
\cdot(1+i)^{2\[n_2/2\]}(1+i)^{2\{n_2/2\}-l_2}
\\ &\qquad
=2^{n_1+\[n_2/2\]-m+l_0}i^{\[n_2/2\]}
(1+i)^{2\{n_2/2\}}(1-i)^{l_2}
\in\mathbb Z[i],
\\ &
(-1)^{l_0+l_1+l_2}A_{\bl}a_2^{l_0+l_1+l_2-m}
=2^{n_1-l_1}(1+i)^{n_2+l_0+l_1-m}
\\ &\qquad
=(1+i)^{n_1-l_1}(1-i)^{n_1-l_1}(1+i)^{n_2+l_0+l_1-m}
\\ &\qquad
=(1+i)^{n_1+n_2-m+l_0}(1-i)^{n_1-l_1}
\in\mathbb Z[i],
\end{align*}
provided that $n_1+\[n_2/2\]-m\ge0$ and $n_1+n_2-m\ge0$
(i.e., that $\alpha_1+\alpha_2/2\ge\alpha$).

\subsection{Measure for $\pi$}
\label{sec:2.3}
Thus, Hata's choice $a_1=2$, $a_2=1+i$ and
$n_0=n_1=n_2=2n$, $m=3n$ with the help of Proposition~\ref{pr:2}
gives the following result.

\begin{theorem}
\label{th:2}
The irrationality exponent of any nonzero
$\gamma\in\mathbb Q\log2+\mathbb Q\pi$
satisfies the inequality
$$
\mu(\gamma)\le8.01604539\dotsc.
$$
\end{theorem}

We would like to refer the interested reader to the notes~\cite{Be}
that could give some feelings of how difficult is evaluating
the irrationality measure of~$\pi$.

\subsection{Double hypergeometric series}
\label{sec:2.4}
Here we present a connection of Hata's construction
with hypergeometric series (that were a major tool
in Section~\ref{sec:1}).

For simplicity, we will set $a=a_1$, $b=a_2$ and deal with
the integrals
$$
J=\int_1^a\frac{(z-1)^{n_0}(z-a)^{n_1}(z-b)^{n_2}}{z^{m+1}}\,\d z
$$
and
$$
J^*=\int_1^b\frac{(z-1)^{n_0}(z-a)^{n_1}(z-b)^{n_2}}{z^{m+1}}\,\d z
$$
giving the simultaneous approximations to $\log a$ and $\log b$.
Applying the starting change of variable $z=1-(1-a)x$ to the first
integral we obtain the single integral
\begin{equation}
J=(-1)^{n_0+1}(1-a)^{n_0+n_1+1}(1-b)^{n_2}
\int_0^1\frac{x^{n_0}(1-x)^{n_1}\Bigl(1-\dfrac{1-a}{1-b}x\Bigr)^{n_2}}
{(1-(1-a)x)^{m+1}}\,\d x
\label{eq:19}
\end{equation}
that may be identified with the Appell hypergeometric function
\begin{align*}
J
&=(-1)^{n_0+1}(1-a)^{n_0+n_1+1}(1-b)^{n_2}
\frac{\Gamma(n_0+1)\,\Gamma(n_1+1)}{\Gamma(n_0+n_1+1)}
\\ &\qquad\times
F_1\biggl(n_0+1;m+1,-n_2;n_0+m+2;1-a,\frac{1-a}{1-b}\biggr)
\end{align*}
(see \cite{Ba}, Section~9.3, formula~(4)),
where the series
$$
F_1(A;B,B';C;X,Y)
=\sum_{\nu=0}^\infty\sum_{\mu=0}^\infty
\frac{(A)_{\nu+\mu}(B)_\nu(B')_\mu}{\nu!\,\mu!(C)_{\nu+\mu}}
X^\nu Y^\mu
$$
is absolutely convergent in the domain $|X|<1$, $|Y|<1$.

The next change of variable
$$
x=(1-y)\bigg/\biggl(1-\frac{1-a}{1-b}y\biggr)
$$
in~\eqref{eq:19} gives the integral representation
\begin{align}
J
&=(-1)^{n_0+m}(1-a)^{n_0+n_1+1}(1-b)^{n_0+n_2+1}(a-b)^{n_1+n_2+1}
\nonumber\\ &\qquad\times
\int_0^1\frac{y^{n_1}(1-y)^{n_0}\,\d y}
{\bigl(a(1-b)-b(1-a)y\bigr)^{m+1}
\bigl((1-b)-(1-a)y\bigr)^{n_0+n_1+n_2-m+1}}
\label{eq:20}
\\
&=(-1)^{n_0+m}\frac{(1-a)^{n_0+n_1+1}(a-b)^{n_1+n_2+1}}
{a^{m+1}(1-b)^{n_1+1}}
\,\frac{\Gamma(n_0+1)\,\Gamma(n_1+1)}{\Gamma(n_0+n_1+1)}
\nonumber\\ &\qquad\times
F_1\biggl(n_1+1;m+1,n_0+n_1+n_2-m+1;n_0+n_1+2;
\frac{b(1-a)}{a(1-b)},\frac{1-a}{1-b}\biggr).
\nonumber
\end{align}
The case $a=2$, $b=1+i$ gives us the following arguments
of the last $F_1$-series:
$$
\frac{1-a}{1-b}=-i=e^{-\pi i/2},
\qquad
\frac{b(1-a)}{a(1-b)}=\frac1{\sqrt2}e^{-\pi i/4}.
$$

Finally, the above changes of variable applied to the integral~$J^*$
produce the same integrals as in~\eqref{eq:19} and \eqref{eq:20}
but with integrations over smooth paths from~$0$ to $(1-b)/(1-a)$
and from~$\infty$ to~$1$, respectively.

\section{Irrationality measure for $\log3$ (after G.~Rhin)}
\label{sec:3}

\subsection{Preliminary remark}
\label{sec:3.1}
As mentioned, the method of Section~\ref{sec:2} have several other applications.
For instance, the choice $a=4/3$, $b=3/2$ and
$n_0=n_1=n_2=2n$, $m=3n$ (cf\. Section~\ref{sec:2.3})
with the help of Proposition~\ref{pr:2} implies that
the irrationality exponent of $\gamma\in\mathbb Q\log2+\mathbb Q\log3$
satisfies the inequality
$\mu(\gamma)\le11.1017577\dots$ (see \cite{Hu}, Corollary~3.1).

\subsection{Back to rational approximations to $\log2$}
\label{sec:3.2}
As we already know from Section~\ref{sec:1.1}, for our starting
integral~\eqref{eq:1} in the case $a=2$ we have
$$
D_n\int_0^1\biggl(\frac{x(1-x)}{1+x}\biggr)^n\frac{\d x}{1+x}
\in\mathbb Z\log2+\mathbb Z,
$$
hence
$$
D_n\int_0^1\biggl(\frac{x(1-x)}{1+x}\biggr)^k\frac{\d x}{1+x}
\in\mathbb Z\log2+\mathbb Z
$$
for any non-negative integer $k\le n$. Considering linear
combinations of the latter integrals we arrive at general inclusions
\begin{equation}
D_n\int_0^1G_n\biggl(\frac{x(1-x)}{1+x}\biggr)\frac{\d x}{1+x}
\in\mathbb Z\log2+\mathbb Z
\label{eq:21}
\end{equation}
valid for all polynomials $G_n(y)\in\mathbb Z[y]$ of degree
$\deg G_n\le n$. To guess a `nice' choice for the polynomial $G_n$,
we start with notifying that
$$
\int_0^1\biggl(\frac{x(1-x)}{1+x}\biggr)^n\frac{\d x}{1+x}
=C\int_0^b\biggl(\frac{x(1-x)}{1+x}\biggr)^n\frac{\d x}{1+x},
$$
where $C$ is a constant (in our case $C=2$) and $b$ is the saddle
point for the integrand: $b=\sqrt2-1$; therefore,
$$
\int_0^1\biggl(\frac{x(1-x)}{1+x}\biggr)^n\frac{\d x}{1+x}
=C\int_0^{(\sqrt2-1)^2}y^nx(y)\,\d y,
$$
where $y=x(1-x)/(1+x)$ and
$x(y)\:(0,(\sqrt2-1)^2)\to(0,b)$ is the inverse function.
Finally,
$$
\int_0^1G_n\biggl(\frac{x(1-x)}{1+x}\biggr)\frac{\d x}{1+x}
=C\int_0^{(\sqrt2-1)^2}G_n(y)x(y)\,\d y;
$$
thus, evaluating the required asymptotic, using inclusions~\eqref{eq:21}
and applying Proposition~\ref{pr:1} result in the estimate
$\mu(\log2)\le1+C_1/C_0$, where
$$
\begin{aligned}
C_0&=-1-\lim_{n\to\infty}
\log\max_{0\le y\le(\sqrt2-1)^2}\bigl\{|G_n(y)|^{1/n}\bigr\},
\\
C_1&=\phantom-1+\lim_{n\to\infty}
\log\max_{0\le y\le(\sqrt2+1)^2}\bigl\{|G_n(y)|^{1/n}\bigr\}.
\end{aligned}
$$
One might now think to look for a polynomial $G_n\in\mathbb Z[y]$
of degree $\le n$ admitting the minimum for the quantity $C_1/C_0$.
Unfortunately, the (non-linear!) problem seems
to be very hard for being solved.

The idea of Rhin \cite{Rh1}, \cite{Rh2}, who introduced
the above construction, was to `linearize' the optimization.
He suggested to look for a polynomial $G^*\in\mathbb Z[y]$
of degree $\le n^*$, say, which is close enough to the optimal
polynomial choice in the problem
\begin{equation}
\min\doublesb{G\in\mathbb Z[y]}{1\le\deg G\le n^*}
\max_{0\le y\le(\sqrt2-1)^2}\bigl\{|G(y)|^{1/n^*}\bigr\},
\label{eq:22}
\end{equation}
and then take $G_n(x)$ to be $(G^*(x))^{\[n/n^*\]}$
for $n$ sufficiently greater than~$n^*$.
For instance, the fact $(\sqrt2-1)^2\approx1/6$ gives one
the first non-trivial approximation $G^*(y)=y^6(6y-1)$
in the problem.

The problem of minimizing the quantity \eqref{eq:22}
is deeply related to evaluating the
$\mathbb Z$-transfinite diameter of the segment $[0,(\sqrt2-1)^2]$.
(The $\mathbb Z$-transfinite diameter of the set $Y\subset\mathbb R$
is defined by the formula
$$
t_{\mathbb Z}(Y)=\inf\doublesb{G\in\mathbb Z[x]}{\deg G\ge1}
\max_{y\in Y}\bigl\{|G(y)|^{1/\deg G}\bigr\},
$$
see \cite{Am1} for problems of computing the quantity.)
This relationship is described in \cite{Am2}; there one can
also find the result $\mu(\log2)<3.991$, which may be achieved
by the method. The latter estimate looks rather close to
the inequality in Theorem~\ref{th:1}; however, it seems to be very
`computer dependent'.

The above method may be used in situations \`a la Section~\ref{sec:2}
as well. For example, we may go back
to simultaneous $\mathbb Z[i]$-approximations
to $\log a_1$ and $\log a_2$ and write
\begin{align*}
\int_1^{a_1}G_n\biggl(\frac{(z-1)^2(z-a_1)^2(z-a_2)^2}{z^3}\biggr)
\,\frac{\d z}z&=B_n\log a_1-B_n',
\\
\int_1^{a_2}G_n\biggl(\frac{(z-1)^2(z-a_1)^2(z-a_2)^2}{z^3}\biggr)
\,\frac{\d z}z&=B_n\log a_2-B_n''
\end{align*}
for any polynomial $G_n(y)\in\mathbb Z[y]$ of degree $\le n$,
where
$$
d^nB_n,\ d^nD_{3n}B_n',\ d^nD_{3n}B_n''\in\mathbb Z[i],
$$
the integer $d>0$ emanates from denominators to the numbers
$a_1,a_2,a_1^{-1},a_2^{-1}$.
(Using the better inclusions achieved by Hata in~\cite{Ha2}
is in this case rather problematic.) Unfortunately, this
way does not look perspective, again due to the fact that
we are required to `linearize' the appeared optimization problem.

\subsection{Another generalization of the integral in~\eqref{eq:1}}
\label{sec:3.3}
On the other hand, we may perform integration in~\eqref{eq:1}
by putting a general polynomial of degree $\le2n$ in
the numerator of the integrand (in place of $x^n(1-x)^n$).
Of course, in this case the polynomial is required
to satisfy some additional conditions.

Let $a=c/d\in\mathbb Q$ with pairwise coprime $c$ and $d>0$,
and let $\Delta$ be a common multiple of the numbers~$c$ and~$d$.
Suppose that a polynomial $H_n(z)\in\mathbb Z[z]$
of degree $\le2n$ may be represented in the form
\begin{equation}
H_n(z)
=\sum_{\nu=0}^nB_\nu\Delta^{n-\nu}z^\nu
+\sum_{\nu=n+1}^{2n}B_\nu z^\nu,
\qquad\text{where}\quad
B_\nu\in\mathbb Z, \; \nu=0,1,\dots,2n.
\label{eq:23}
\end{equation}
(Clearly, for $a=2$ the polynomial $H_n(z)=(z-1)^n(z-2)^n$
has the desired form.) Then for the integral
$$
I(n)=(1-a)\int_0^1\frac{H_n(d-d(1-a)x)}{d^n(1-(1-a)x)^{n+1}}\,\d x
$$
we deduce
\begin{align*}
I(n)
&=\sum_{\nu=0}^nB_\nu\Delta^{n-\nu}d^{\nu-n}
(1-a)\int_0^1(1-(1-a)x)^{\nu-n-1}\,\d x
\\ &\qquad
+\sum_{\nu=n+1}^{2n}B_\nu d^{\nu-n}
(1-a)\int_0^1(1-(1-a)x)^{\nu-n-1}\,\d x
\\
&=\sum_{\nu=0}^{n-1}B_\nu\Delta^{n-\nu}d^{\nu-n}\frac{a^{\nu-n}-1}{n-\nu}
-B_n\log a
-\sum_{\nu=n+1}^{2n}B_\nu d^{\nu-n}\frac{1-a^{\nu-n}}{\nu-n},
\end{align*}
hence
$$
I(n)\cdot D_n\in\mathbb Z\log a+\mathbb Z.
$$

In general, having a set of $k$ rational numbers
$a_j=c_j/d$ for $j=1,\dots,k$,
we suppose that the polynomial $H_n(z)\in\mathbb Z[z]$
of degree $\le2n$ has representation~\eqref{eq:23}
with $\Delta$ being a multiple of the numbers
$c_1,\dots,c_k,d$. Then setting
\begin{equation}
I(n;a_j)=(1-a_j)\int_0^1\frac{H_n(d-d(1-a_j)x)}{d^n(1-(1-a_j)x)^{n+1}}\,\d x,
\qquad j=1,\dots,k,
\label{eq:24}
\end{equation}
we obtain
$$
I(n;a_j)\cdot D_n
=-B_n\log a_j+A_{nj}\in\mathbb Z\log a_j+\mathbb Z,
\qquad j=1,\dots,k,
$$
again simultaneous approximations to $\log a_1,\dots,\log a_k$.
(In fact, the choice
$$
H_n(z)=\Delta^{2n}(z-1)^{\[\beta_0n\]}(z-a_1)^{\[\beta_1n\]}
\dotsb(z-a_k)^{\[\beta_kn\]}
$$
where $\beta_j=\alpha_j/\alpha$ for $j=0,1,\dots,k$,
gives us exactly the same approximations as in Section~\ref{sec:2}.
The case $\beta_1=\dots=\beta_k$ was previously treated
in~\cite{Rh1} and~\cite{RT}.)

Finding a suitable polynomial $H_n(z)$ for a given set of
the numbers $a_1,\dots,a_k$ is very similar to that of Section~\ref{sec:3.2}.
The change of variable $z_j=d-d(1-a_j)x$ in the integrals~\eqref{eq:24}
(hence, integrating
then a simpler expression over the segment $[d,da_j]$)
leads to the problem of finding a polynomial $H_n(z)\in\mathbb Z[z]$
of degree $\le2n$ with expansion~\eqref{eq:23} such that
the quantity
$$
\max_{z\in Z}\biggl\{\biggl|\frac{H_n(z)}z\biggr|^{1/n}\biggr\},
\qquad Z=\bigcup_{j=1}^k[d,da_j],
$$
is as small as possible. The algorithmic solution to this
optimization problem by means of the LLL-algorithm was
recently proposed by Q.~Wu~\cite{Wu}. This gives one a machinery
to produce fairly good estimates for linear forms
in the logarithms of {\it rational\/} numbers.

\subsection{Measure for $\log3$}
\label{sec:3.4}
To derive a nice irrationality measure for $\log3$,
Rhin constructs in~\cite{Rh2} simultaneous approximations
to the logarithms of $a_1=2/3$, $a_2=4/3$
and use the following (very complicated)
choice of the polynomial~\eqref{eq:23}:
\begin{align*}
&
H_n(z)
=2^{14}\cdot3^{2n+7}\cdot(z-1)^{\[0.704324n\]}
\Bigl(z-\frac23\Bigr)^{\[0.552418n\]}
\Bigl(z-\frac43\Bigr)^{\[0.447582n\]}
\\ &\qquad\times
(5z-4)^{\[0.109072n\]}
(17z^2-34z+16)^{\[0.038934n\]}
(19z^2-36z+16)^{\[0.054368n\]}
\end{align*}
(a `justification' of the choice is done in~\cite{Wu}).
By these means he proves

\begin{theorem}
\label{th:3}
The irrationality exponent of any nonzero
$\gamma\in\mathbb Q\log2+\mathbb Q\log3$
satisfies the inequality
$$
\mu(\gamma)<8.616.
$$
\end{theorem}

Further results in this direction (e.g., irrationality
measures for $\log5$, $\log7$ etc.)
may be found in~\cite{Wu}.

\section{Concluding improvisations}
\label{sec:4}

Connections with the hypergeometric subject (indicated
in Sections~\ref{sec:1.1} and \ref{sec:2.4} above)
could play a role in further
improvements of the irrationality measures of logarithms
and related constants. For instance, Euler's transform
(see, e.g., \cite{Ba}, Section~2.4, formula~(1))
$$
{}_2F_1\biggl(\begin{matrix} A, \, B \\ C \end{matrix}
\biggm|z\biggr)
=\frac1{(1-z)^A}\cdot{}_2F_1\biggl(\begin{matrix} A, \, C-B \\ C \end{matrix}
\biggm|\frac{-z}{1-z}\biggr)
$$
translates the value $z=1-a=-1$ of Section~\ref{sec:1} into
$-z/(1-z)=1/2$. This leads to a ${}_2F_1$-series
with {\it positive\/} terms and makes possible
the analytic evaluation of the quantity~\eqref{eq:3}
without using the integral representation~\eqref{eq:2}---we
may get rid of the integral (the idea belongs to K.~Ball,
cf\. \cite{Zu}, the proof of Lemma~4). However, other hypergeometric
ingredients are required for real improvements.

We find quite curious that Ramanujan's formulae for~$\pi$,
in particular
\begin{equation}
\begin{gathered}
\sum_{\nu=0}^\infty\frac{(1/4)_\nu(1/2)_\nu(3/4)_\nu}{\nu!^3}
(21460\nu+1123)\cdot\frac{(-1)^\nu}{882^{2\nu+1}}
=\frac4\pi,
\\
\sum_{\nu=0}^\infty\frac{(1/4)_\nu(1/2)_\nu(3/4)_\nu}{\nu!^3}
(26390\nu+1103)\cdot\frac1{99^{4\nu+2}}
=\frac1{2\pi\sqrt2}
\end{gathered}
\label{eq:25}
\end{equation}
(see \cite{Ra}, equations (39) and (44)) and several others,
might be used for constructing good
rational approximations to $\pi$ and $\pi\sqrt{d}$, where
$d$~is a positive integer. Namely, one can expect reasonable
estimates for the corresponding irrationality measures
by constructing explicit Pad\'e approximations (of either
first or second type) to the functional system
$1$, $f(z)$, $f'(z)$, $f''(z)$, where
$$
f(z)
=\sum_{\nu=0}^\infty\frac{(1/4)_\nu(1/2)_\nu(3/4)_\nu}{\nu!^3}z^\nu
={}_3F_2\biggl(\begin{matrix} \frac14, \, \frac12, \, \frac34 \\
1, \, 1 \end{matrix}
\biggm|z\biggr).
$$
The paper \cite{Ne} provides Pad\'e approximations to the
homogeneous system $f(z)$, $f'(z)$, $f''(z)$ (without~$1$)
that are not enough for our purposes. Finally, we should
mention that a general result of A.~Galochkin in~\cite{Ga}
(proved by a proper variation of Siegel's method) yields
the qualitative linear independence of the numbers
$1$, $f(1/b)$, $f'(1/b)$, and $f''(1/b)$ for integers~$b$
satisfying $|b|>b_0$, where the value of~$b_0$ is so huge
that $b=-882^2$ and $b=99^4$ in~\eqref{eq:25} do not suit.

\begin{acknowledgements}
I thank G.~Rhin kindly for introducing me to the subject
of transfinite diameters and their number-theoretic applications,
in particular those presented in Section~\ref{sec:3}. Special
gratitude is due to P.~Bundschuh, the fruitful discussions
with whom during my long-term stay at Cologne University
were crucial for this writing. I thank J.~Guillera
for attracting my attention to Ramanujan-type formulae and
making me familiar with the manuscript~\cite{Gu},
and J.~Sondow for several suggestions.
\end{acknowledgements}


\noindent
\hbox to70mm{\vbox{\hsize=70mm%
\leftline{Moscow Lomonosov State University}
\leftline{Department of Mechanics and Mathematics}
\leftline{Vorobiovy Gory, GSP-2}
\leftline{119992 Moscow, RUSSIA}
\leftline{{\it E-mail\/}: \texttt{wadim@ips.ras.ru}}
\leftline{{\it URL\/}: \texttt{http://wain.mi.ras.ru/}}
}}

\end{document}